\documentclass[a4paper]{amsart}

\usepackage{amssymb}
\usepackage{amsmath,amsthm}
\usepackage{enumerate,paralist}
\usepackage{amstext}
\usepackage{psfrag}
\usepackage{dsfont}
\usepackage[dvips]{epsfig}
\usepackage[english]{babel}
\usepackage{cite}
\usepackage{hyperref}
\usepackage{color}
\usepackage{cleveref}
\usepackage[svgnames]{xcolor}
\hypersetup{hidelinks,colorlinks=true,allcolors=DarkBlue}

\usepackage{cite}

\theoremstyle{plain}
\newtheorem{theorem}{Theorem}

\newtheorem{lemma}{Lemma}
\theoremstyle{definition}

\def\diag{{\rm diag}}
\def\sign{\,{\rm sign}}

\begin{document}
\frenchspacing

\title[Asymptotically optimal feedback control]
{Asymptotically optimal feedback control \\ for a system of linear oscillators}

\author{Alexander Ovseevich}
\address
{
Institute for Problems in Mechanics, Russian Academy of Sciences \\
119526, Vernadsky av., 101/1, Moscow, Russia.
} 
\email{ovseev@ipmnet.ru}

\author{Aleksey Fedorov}
\address
{
Institute for Problems in Mechanics, Russian Academy of Sciences \\
119526, Vernadsky av., 101/1, Moscow, Russia.
}
\email{akfedorov@student.bmstu.ru}

\maketitle

\begin{abstract}
Problem of damping of an arbitrary number of linear oscillators under common bounded control is considered. 
We are looking for a feedback control steering the system to the equilibrium. 
The obtained control is asymptotically optimal: the ratio of motion time to zero with this control to the minimum one is close to unity, 
if the initial energy of the system is large. 
Some of the results based on the new lemma about sustainable observability autonomous linear systems.


\end{abstract}

\section{Introduction} 

One of the classical achievements in control theory is the explicit construction of the minimum time damping for a single linear oscillator. 
In this paper, we consider the next in complexity problem of damping of an arbitrary number of oscillators under a common bounded control. 
Probably, in this case an explicit construction of the optimal feedback control is impossible, and even the numerical solution is a hard problem. 
We are looking for a non-optimal feedback control which brings the system to the equilibrium point. 
The control obtained is asymptotically optimal: the ratio of the duration of steering to zero under our control to the minimum one is close to unity, 
if the initial energy is large.

\section{The problem statement} 
Equation of motion for the control system of $N$ linear oscillators with eigenfrequencies $\omega_i$ has the following form:
\begin{eqnarray}\label{syst1}
	&\dot{x}={A}x+{B}u,
	\quad 
	x=(x_1,y_1,\dots,x_N,y_N)^*\in{\mathbb{R}}^{2N},
	\quad 
	u\in {\mathbb{R}},
	\quad
	|u|\leq1,\\[1em]
\label{syst2} &{A} = \left( {\begin{array}{*{20}c}
	0	&	1	&	&	&	\\
	{ - \omega_1^{2}} &	0	&	&	&	\\
	&	& \ddots	&	&	\\
	&	&	&	0	&	1	\\
	&	&	&	{-\omega_N^{2}}	&	0	\\
\end{array} } \right), 
\qquad
	{B}=\left( \begin{array}{c}
	0  \\
	1  \\
	\vdots  \\
	0  \\
	1  \\
\end{array}\right).
\end{eqnarray}

Solution of the linear minimum time problem is completely equivalent to the boundary value problem for the Pontryagin maximum principle corresponding to the Hamiltonian
$h(x,\psi)=(Ax,\psi)+|B^*\psi|-1=\max \{(Ax,\psi)+(Bu,\psi)-1\}$, 
where $\psi\in\mathbb{R}^{2N}$ is the vector of adjoint variable (momentum), $u=\sign(B^*\psi)$. 
In particular, knowledge of the momenta defines the control uniquely.

Geometrically, the maximum principle says that the momentum (the vector of adjoint variables) $\psi$ is the inner normal to the reachable set $\mathcal{D}(T(x))$ 
($T(x)$ is the time of steering $x$ from zero point). 
We would like to use the inner normal to {\it approximation} of the reachable set as the momentum.

There are other possible methods of feedback control design, based, e.g., upon the Kalman approach to the feedforward control for linear systems \cite{chern,ovseev3}.

\section{The proposed method} 
Three strategies are used coherently in our method. 
At high energy zone, normal to an approximation of the reachable set, which is close to the real one in the long run \cite{ovseev, ovseev2}, is used as a momentum. 
The control obtained  can be applied within the small energy zone as well, but then its quasioptimal properties are lost. 
Moreover, the control affect the system like a dry friction, so that in some states it does not allow moving at all. 
There are another scenarios of more general character: the motion might occur in a vicinity of a limit set (attractor) not containing the equilibrium state. 
By applying the control with a reduced upper limit (second type of control) allows to delay undesirable pulling in the attractor. 
This allows the system to reach a sufficiently small vicinity of zero point, where the final stage of the control should be applied.

At the terminal, third stage an approach to local feedback control of \cite{anan, korobov}, based on common Lyapunov functions, is used. 
This method works in a sufficiently small vicinity of zero. 
In order to get to this immediate vicinity, it is necessary that it contains the entire zone of pulling into  attractor for the preceding type of  control. 
The reduction of the zone achieved at the second stage is sufficient for this purpose.

\section{Control in high energy zone}\label{high}
One of the main results of \cite{ovseev}, applicable to the system of $N$ oscillators is this: 
the reachable set $D(T)$ equals asymptotically as $T\to\infty$ to the set  $T\Omega,$ where $\Omega$ is a fixed convex body. 
More precisely, suppose that the momentum $p$ is written in the form $p=(p_i)$, 
where $p_i=(\xi_i,\eta_i),$ ${i=1,\dots,N}$, $\xi_i$ is the dual variable for $x_i$, $\eta_i$ is the dual variable for $y_i$, and $z_i=(\eta_i^2+{\omega_i^{-2}}{\xi_i^{2}})^{1/2}$. 
In the non-resonant case (there are no nontrivial relations  $\sum{m_i\omega_i}=0,\, m_i\in\mathbb{Z}$), 
the support function $H_T$ of the reachable set $D(T)$ has as $T\to\infty$ the asymptotic form:
\begin{equation}\label{approxN0}
	{H}_T(p)=\frac{T}{(2\pi)^N}\int_{0}^{2\pi}\dots\int_{0}^{2\pi}\left|\sum_{i=1}^Nz_{i}\cos\varphi_{i}\right|d\varphi_{1}\dots d\varphi_{N}+o(T)=T\mathfrak{H}(z)+o(T),
\end{equation}
the support function of the compact  $\Omega$ is given by the main term $\mathfrak{H}(z)$. 
If $N=1$ we get $\mathfrak{H}(z)=\frac2\pi|z|$, 
if $N=2$ the function $\mathfrak{H}$ can be expressed via elliptic integrals.

The basic idea of our control design is to substitute the set $T\Omega$  for $D(T)$, 
and the normal to this set for  momentum. 
If the phase vector $x$ belongs to the boundary of $T\Omega$, then
\begin{equation}\label{approx3}
	T^{-1}x=\frac{\partial {H}_\Omega}{\partial p}(p)
\end{equation} 
for a momentum $p=p(x)$. 
Note that the support function ${H}_{\Omega}$ is differentiable, and Eq. (\ref{approx3}) has a unique solution, 
because of the smoothness of the boundary of $\Omega$ \cite{ovseev2}. 
Our feedback control is given by
\begin{equation}\label{approx4}
	u(x)=-\sign(B,p(x)).
\end{equation}

\subsection{Asymptotic optimality of control (\ref{approx4})}\label{quasioptimal} 

We define a polar-like coordinate system (if $N=1$ we get the canonical polar coordinate system in a plane),
well suited for repre\-sen\-ta\-tion of  motion under the control $u$. 
Write the phase vector $x$ in the form $x=\rho\phi$, where $\rho>0$ and $\phi\in\omega=\partial\Omega$. 
In terms of Eq. (\ref{approx4}), $\rho=T$ and $\phi=\frac{\partial {H}_\Omega(p)}{\partial p}$. 
In these coordinates equations of the motion have the form
\begin{equation}\label{T}
	\dot \rho=-\left|\left(\frac{\partial {\rho}}{\partial x},B\right)\right|,
	\quad
	\dot\phi=A\phi+\frac{1}{\rho}\left(Bu+\phi \left|\left(\frac{\partial {\rho}}{\partial x},B\right)\right|\right).
\end{equation}
An eikonal-type equation holds for the function $\rho=\rho(x)$
\begin{equation}\label{Euler}
	H_\Omega\left(p\right)=1,\quad p=\frac{\partial \rho}{\partial x}.
\end{equation}
It is ``dual'' to the equation $\rho({\partial {H_\Omega}}/{\partial p})=1$ of the surface $\omega$.
Eq. (\ref{Euler}) can be used for averaging the right-hand side of the first equality (\ref{T}) in time,
and is the basis of the proof of the following statement on the asymptotic optimality of control (\ref{approx4}):
\begin{theorem}\label{main_approx}
	Consider evolution (\ref{T}) of  $\rho$ under control (\ref{approx4}). 
	Put $M=\min\{\rho(0),\rho(T),T\}$. 
	Then as $M\to+\infty$ we have
	\begin{equation}\label{approx_T}
		{(\rho(0)-\rho(T))}/{T}=1+o(1).
	\end{equation} 
	Under any other admissible control
	\begin{equation}\label{approx_T2}
		{(\rho(0)-\rho(T))}/{T}\leq 1+o(1).
	\end{equation}
\end{theorem}

\subsection{Comparison with the maximum principle}\label{maximum} 

One can approach the issue of asymptotic optimality of control (\ref{approx4}) by comparison of the differential equations of the motion 
under the control with equations of the Pontryagin maximum principle. 
This requires understanding the time-evolution of the momentum $p(x)$, involved in (\ref{approx4}). 
A description is given by the following equation:
\begin{equation}\label{attractor_syst22p1}
	\dot p=-A^*p+\widetilde Bu, \mbox{ where } \widetilde B=\frac{\partial^2\rho}{\partial x^2}B.
\end{equation}
note that if the latter equation would not contain the second term $\widetilde Bu$, 
then the equation for  $\psi=-p$ would coincide with with the maximum principle equation for adjoint variables. 
However, the matrix $\frac{\partial^2\rho}{\partial x^2}$ is a homogeneous function of $x$ of degree $-1$, 
therefore the said second term has order $O(\frac1{|x|})$ for $x$ large, and, therefore, is small. 
Note that the maximum condition $u=\sign(B,\psi)=-\sign(B,p)$ holds for control (\ref{approx4}). 
It remains to find out to what extent the condition $h(x,\psi)=(Ax,\psi)+|B^*\psi|-1=0$ holds. 
In fact, the arguments of the previous section imply that $h(x,\psi)=0$ ``on the average'' in the non-resonant case. 
Indeed, $(Ax,\psi)=-(Ax,{\partial\rho}/{\partial x})=0,$ and the average value of $|B^*\psi|=|B^*p|$ is close to 1 for $x$ sufficiently large, according to Theorem \ref{main_approx}.

Thus, the maximum principle equation for the vector $(x,\psi)$, 
where $\psi=-\frac{\partial\rho}{\partial x}$, holds ``on the average'' with a small error as $x$ is large.

\subsection{Efficiency of control (\ref{approx4}) in the near zone}\label{efficiency} 

In accordance with Theorem \ref{main_approx} the time of motion from the level set $\rho=M$ to the level set $\rho=N$ under control (\ref{approx4}) is asymptotically $(M-N)(1+o(1))$, 
if $M,N,$ and $M-N$ are very large. 
Now we show that a nonasymptotical estimate holds: the time of motion $T$ is $O(M-N)$, if $M,N$ and $M-N$ are greater than a constant $C(A,B)$, 
depending only on parameters of our system of oscillators. 
Relation (\ref{T}) reduces the required estimate to the inequality
\begin{equation}\label{T33}
	\int_0^T|(p,B)|dt\geq cT,
\end{equation}
where  $c=c(A,B)$ is a positive constant. 
To prove (\ref{T33}), we use the following important lemma on completely controllable time-invariant linear systems.
\begin{lemma}\label{observation2}
	Suppose that $\dot x=\alpha x,\,y=\beta x$ is a completely controllable time-invariant linear system. 
	Then, a solution  $z$ of equation $\dot z=\alpha z+f$ in the interval $I$ of an {\bf integer} length $T\geq1$ satisfies a priori estimate 
	$\int_I|z|dt\ll \int_I|\beta z|dt+\int_I|f|dt$ (here, $\ll$ is the Vinogradov symbol).
\end{lemma}
We apply the Lemma to Eq. (\ref{attractor_syst22p1}) with the phase vector $p$, observation $y=(p,B)$, and the right-hand side $f={\widetilde B}u$. 
Assume, that in the entire time interval $I$ of integer length $T$ the motion governed by $\dot{x}={A}x+{B}u$ takes place within the domain $\rho(x)\geq C$. 
Then $|f|=O(1/C)$ in the entire interval. 
In addition, Eq. (\ref{Euler}) holds for $p$ and therefore $T\ll\int_I|p|dt$. 
The estimate of Lemma \label{observation2} and the eikonal Eq. (\ref{Euler}) gives that $T\ll\int_I|p|dt\ll \int_I|(p,B)|dt+\frac1C T$. 
By taking a sufficiently large constant $C=C(A,B)$, we obtain that $T\ll\int_I|(p,B)|dt$, and this is inequality (\ref{T33}) in another notation. 
Thus, we have:
\begin{theorem}\label{observation5}
	Suppose that the motion from the level set $\rho=M$ to the level set $\rho=N$ under  control (\ref{approx4}) goes within the domain $\rho(x)\geq C(A,B)$, 
	in the time interval of integer length $T$, where $C(A,B)$ is a (sufficiently large) constant, depending on our system of oscillators only.
	Then $T\leq c(M-N)$, where $c=c(A,B)$ is a strictly positive constant.
\end{theorem}
For the reduced control
\begin{equation}\label{controlU}
	u_U(x)=Uu(x),\,|U|\leq1.
\end{equation}
we obtain the following result.
\begin{theorem}\label{observation52}
	Suppose that the motion from the level set $\rho=M$ to the level set $\rho=N$ under control (\ref{approx4}) goes within the domain $\rho(x)\geq UC(A,B)$, 
	in the time integer of integer length $T$, where $C(A,B)$ is a (sufficiently large) constant from Theorem~\ref{observation5}. 
	Then, $T\leq\frac{c}{U}(M-N)$, where $c=c(A,B)$ is another constant from Theorem \ref{observation5}.
\end{theorem}

\subsection{Singular motion}\label{attractor} 

According to Eqs.~(\ref{T}) the value of $\rho$ under control (\ref{approx4}) does not increase, 
but might stay constant if the condition $({\partial\rho}/{\partial x},B)=(p,B)=0$ holds in a time interval. 
In particular, this condition is fulfilled along any $\omega$-limit set (attractor) not containing the equilibrium.

Consider the ``dual'' dynamic system describing the motion of the vector $p=\frac{\partial\rho}{\partial x}(\phi)$. 
Put $\widetilde B=\frac{\partial^2\rho}{\partial x^2}B$.
We obtain, according to formula (\ref{attractor_syst22p1}), that
\begin{equation}\label{attractor_syst22p}
	\dot p=-A^*p+\widetilde Bu, \quad (p,B)=0,
\end{equation}
from which we can get the value $u=(p,AB)/(\widetilde B,B)$ for the singular control. 
The expression $\widetilde B=\frac{\partial^2\rho}{\partial x^2}B$ can also be written as a function of momentum $p$: 
$\widetilde{B}=\left(\frac{\partial^2 {H}}{\partial p^2}\right)^{-1}B$. 
Therefore, the motion along the attractor is described by the dynamical system
\begin{equation}\label{p_attractor}
	\dot p=-A^*p+{\widetilde B}f(p), \quad (p,B)=0,\mbox{ where }f(p)=(p,AB)/(\widetilde B,B)
\end{equation}
on the topological sphere $\widetilde \sigma=\left\{H(p)=1,\left(p,B\right)=0\right\}$ of dimension $2N-2$. 
Namely, the phase trajectory of system (\ref{p_attractor}) is contained in an attractor, iff the inequality $|f(p)|\leq1$ holds within the trajectory.
\begin{theorem}
	The number $\mu$, defined as the minimum over trajectories of  \ref{p_attractor} of the maximum of the function $|f|$ on a trajectory, is strictly positive.
\end{theorem}
Importance of $\mu$ is due to the fact that it gives an exact bound for the efficiency zone of control (\ref{approx4}):
\begin{theorem}\label{mu} 
	Suppose that $\epsilon>0$, and the motion under control (\ref{approx4}) in a sufficiently long time interval $[a,b]$ of length $T$ goes within the domain $\rho\geq\mu^{-1}+\epsilon$. 
	Then $\rho(a)-\rho(b)\geq c(\epsilon)T$, 
	where $c(\epsilon)$ is a positive constant. 
	On the other hand, there are infinitely long motions within the domain $\mu^{-1}-\epsilon\leq\rho\leq\mu^{-1}$, 
	where $\rho$ stays constant.
\end{theorem}
In notations of Theorem \ref{observation5} this means that $C(A,B)=\mu^{-1}+\epsilon$.

\section{The feedback nearby the terminal point}\label{terminal} 

The design of our local feedback control goes back to \cite{korobov}, and it uses a preliminary reduction of the system (\ref{syst1})--(\ref{syst2}) to a canonical form by means of transformations
\begin{equation}\label{transformations}
	A\mapsto A+BC,\quad  u\mapsto u-Cx ,\quad A\mapsto D^{-1}AD,\quad B\mapsto D^{-1}B,
\end{equation}
corresponding to adding a linear feedback control, and to coordinate changes (gauge transformations). 
We state the result as follows:
\begin{lemma}\label{canonical} 
	By transformations (\ref{transformations}), the system (\ref{syst1})--(\ref{syst2}) can be reduced to the form 
	$\dot {\mathfrak x}={\mathfrak A}{\mathfrak x}+{\mathfrak B}{\mathfrak u}$, where
	\begin{equation}\label{AB}
		\begin{array}{c}
			\mathfrak{A} = 
			\left( {\begin{array}{*{20}c}
			0 & &  &   \\
			-1 & 0 &  &   \\
			& -2 &0&\\
			& &\ddots & \ddots    \\
			&&  & -2N+1 & 0  \\
			\end{array} } \right), \quad
			\mathfrak{B} = 
			\left( \begin{array}{c}
			1  \\
			0  \\
			0  \\
			\vdots  \\
			0  \\
			\end{array}  \right).
		\end{array}
\end{equation}
To do this, the matrix of the linear feedback should have the form
\begin{equation}\label{C}
	C=(c_1\,0\,c_2\,0\,\dots\, c_N\,0),\,c_k=(-1)^{N+1}\omega_k^{2N}\prod_{i\neq k}(\omega_i^{2}-\omega_k^{2})^{-1}.
\end{equation}
The gauge matrix $D$ transforms the standard basis $e_i=(\delta_{ij})$ of $\mathbb{R}^{2N}$ into the basis
\begin{equation}\label{e}
	\mathfrak{e}_i=\frac{(-1)^{i-1}}{(i-1)!}(A+BC)^{i-1}B,\,i=1,\dots,2N,
\end{equation}
and has the following form: 
Define $2\times2$ matrices
\begin{equation}\label{d}
	d_{ij}=(-1)^{j-1}\lambda_i^{j-1}\left(%
	\begin{array}{cc}
		0 &  -\frac{1}{(2j-1)!} \\
		\frac{1}{(2(j-1))!} & 0
	\end{array}\right), 
	\mbox{ where }
	\lambda_k=\sum_{i\neq k}\omega_i^2.
\end{equation} 
Then
\begin{equation}\label{D}
	D\mbox{ is the $N\times N$ matrix $(d_{ij})$ of $2\times2$ blocks $d_{ij}.$}\end{equation}
\end{lemma}
When regarded as an existence theorem of a canonical form, 
without explicit formulas for matrices $C$ and $D$, 
Lemma \ref{canonical} is a particular case of the Brunovsky theorem \cite{brun}. 
By following \cite{anan}, introduce a matrix function of time, related to system (\ref{AB}):
\begin{equation}\label{delta}
	\delta(T)=\diag(T^1,T^{2}, \dots, T^{2N})^{-1}.
\end{equation}
In what follows the parameter $T$ will be a function $T=T(\mathfrak{x})$ of the phase vector.
Define, in accordance with \cite{korobov,anan}, the matrices
\begin{equation}\label{M}
	\begin{array}{l}
		\mathfrak{q}=(\mathfrak{q}_{ij}),
		\quad 
		\mathfrak{q}_{ij}=\int_0^1 x^{i+j-2}(1-x)dx=[(i+j)(i+j-1)]^{-1}, \\[1em]
		\mathfrak{Q}=\mathfrak{q}^{-1},
		\quad
		\mathfrak{C}=-\frac{1}{2}\mathfrak{B}^{*}\mathfrak{Q}, 
		\quad 
		\mathfrak{M}=\diag(1, 2, \dots, 2N) \\
\end{array}
\end{equation}
and the feedback control by
\begin{equation}\label{u}
	{\mathfrak u}(\mathfrak{x})=\mathfrak{C}\delta(T( \mathfrak{x}))\mathfrak{x},
\end{equation}
where the function $T=T(\mathfrak{x})$ is given implicitly by
\begin{equation}\label{condu}
	(\mathfrak{Q}\delta(T)\mathfrak{x},\delta(T)\mathfrak{x})=\kappa^2=\frac{1}{2N(2N+1).}
\end{equation}
The basic result on steering the canonical system (\ref{AB}) to zero is as follows:
\begin{theorem}\label{main}
	\begin{itemize}
		\item[A)] The matrix $\mathfrak{Q}$ defines a common quadratic Lyapunov's function for the matrices $-\mathfrak{M}$ and $\mathfrak {A+BC}$.
		\item[B)] Eq. (\ref{condu}) defines $T=T(\mathfrak{x})$ uniquely.
		\item[C)] Control (\ref{u}) is bounded: $|{\mathfrak u}|\leq\frac{\kappa}{2}\sqrt{\mathfrak{Q}_{11}} $.
		\item[D)] Control (\ref{u}) brings the point $\mathfrak{x}$ to $0$ in time $T(\mathfrak{x})$.
		\item[E)] The matrix $\mathfrak{Q}$ is an even integer one.
	\end{itemize}
\end{theorem}
The assertions (A)--(D) of Theorem \ref{main} were obtained in \cite{korobov} in a less precise form.

\section{The control matching.}\label{match}
In Sec.~\ref{terminal}, we designed a \textit{local} feedback control, 
which works in a vicinity of zero. 
The switching to this control should happen at the boundary of an \textit{invariant} domain with respect to the phase flow such that the local feedback control can be applied within the interior. 
We confine ourselves with the invariant domains of the form
\begin{equation}\label{theta}
	G_\Theta=\{{\mathfrak x}: T({\mathfrak x})\leq\Theta\}=\{{\mathfrak x}: (\mathfrak{Q}\delta(\Theta){\mathfrak x},\delta(\Theta){\mathfrak x})\leq1\}.
\end{equation}
The invariant domain  $G_\Theta$ should satisfy two conditions:
\begin{itemize}\label{conditionsG}
	\item[A:] The domain $G_\Theta$ contain the inefficiency domain $\{\rho(x)\leq UC(\underline\omega)\}$ of the preceding control.
	\item[B:] The domain  $G_\Theta$ is contained in the strip $\{|Cx|\leq1/2\}$, where $C$ is the matrix (\ref{C}).
\end{itemize}
The condition B allows to use at the terminal stage controls ${\mathfrak u}$ which are less than 1/2 in absolute value. 
Therefore, the constant $\kappa^2$ in (\ref{condu}) should be equal to $\frac1{2N(2N+1)}$. 
If we applied at the preceding stage control (\ref{controlU}), 
the condition A says that the set $UC(A,B)\Omega$ is contained in $G_\Theta$. 
Here, $C(A,B)$ is the estimate found in Section \ref{efficiency} for the ``radius'' of the attractor free domain. 
In other words, the inequality should be fulfilled for the support functions
\begin{equation}\label{condA}
	UC(A,B)H_\Omega(D^*p)\leq(\delta(\Theta)^{-1}\mathfrak{q}\delta(\Theta)^{-1}p,p)^{1/2},
\end{equation}
where $D$ is matrix (\ref{D}). 
It is clear that the inequality holds, provided that $U$ is sufficiently small.

The condition B says, 
that the value of the support function of ellipsoid $G_\Theta$ at the vector ${D^*}^{-1}C$ does not exceed $1/2$ in absolute value. 
In other words,
\begin{equation}\label{condB}
	(\delta(\Theta)^{-1}\mathfrak{q}\delta(\Theta)^{-1}{D^*}^{-1}C,{D^*}^{-1}C)^{1/2}\leq1/2.
\end{equation}
Certainly, the inequality holds for sufficiently small $\Theta$. 
After choice of $\Theta$ we have to choose the bound $U$ for the control at the second stage in accordance with inequality (\ref{condA}). 
Then, condition A and B are met. 
The switching to the third, terminal stage should happen upon arriving at the boundary 
$\{(\mathfrak{Q}\delta(\Theta){\mathfrak x},\delta(\Theta){\mathfrak x})=1\}$ of $G_\Theta$.

\section{The final asymptotic result}\label{asymp}

\begin{theorem}\label{main_approx22}
	Let $T=T(x)$ be the motion time from the initial point $x$ to the equilibrium under our three-stage control, and let $\tau=\tau(x)$ be the minimum time. 
	Then, as $\rho(x)\to+\infty$ we have asymptotic equalities
	\begin{equation}\label{approx_T12}
		\rho(x)/T(x)=1+o(1), \, \tau(x)/T(x)=1+o(1).
	\end{equation}
	\end{theorem}
Note that the equality $\rho(x)/T(x)=1+o(1)$ follows from the asymptotic theory of reachable sets
\cite{ovseev}.

Thus, we designed an asymptotically optimal and, at the same time, constructive feedback control for arbitrarily large system of oscillators.  
The design  reduces to solution of Eq. (\ref{approx3}), 
which is essentially the same as the well-studied problem of maximization of a linear form on a convex hypersurface in $\mathbb{R}^N $.

\section*{Acknowledgements}
The authors thank the anonymous referee for  careful reading of the manuscript and useful comments.
This work was supported by RFBR, project 11-08-00435.

\end{document}